\documentclass[11pt]{amsart}
\usepackage{amsmath,amssymb,mathrsfs}
\usepackage{amsmath,amssymb}
\newtheorem{theorem}{Theorem}[section]
\newtheorem{proposition}[theorem]{Proposition}
\newtheorem{corollary}[theorem]{Corollary}
\newtheorem{lemma}[theorem]{Lemma}
\theoremstyle{definition}
\newtheorem*{definition}{Definition}

\begin{document}

\title[Ricci solitons in higher dimensions]{Rotational symmetry of Ricci solitons in higher dimensions}
\author{Simon Brendle}
\begin{abstract}
Let $(M,g)$ be a steady gradient Ricci soliton of dimension $n \geq 4$ which has positive sectional curvature and is asymptotically cylindrical. Under these assumptions, we show that $(M,g)$ is rotationally symmetric. In particular, our results apply to steady gradient Ricci solitons in dimension $4$ which are $\kappa$-noncollapsed and have positive isotropic curvature.
\end{abstract}
\address{Department of Mathematics \\ Stanford University \\ Stanford, CA 94305}
\thanks{The author was supported in part by the National Science Foundation under grants DMS-0905628 and DMS-1201924.}
\maketitle

\section{Introduction}

This is a sequel to our earlier paper \cite{Brendle2}, in which we proved a uniqueness theorem for the three-dimensional Bryant soliton. Recall that the Bryant soliton is the unique steady gradient Ricci soliton in dimension $3$, which is rotationally symmetric (cf. \cite{Bryant}). In \cite{Brendle2}, it was shown that the three-dimensional Bryant soliton is unique in the class of $\kappa$-noncollapsed steady gradient Ricci solitons:

\begin{theorem}[S.~Brendle \cite{Brendle2}]
\label{dim.3}
Let $(M,g)$ be a three-dimensional complete steady gradient Ricci soliton which is non-flat and $\kappa$-noncollapsed. Then $(M,g)$ is rotationally symmetric, and is therefore isometric to the Bryant soliton up to scaling. 
\end{theorem}

Theorem \ref{dim.3} resolves a problem mentioned in Perelman's first paper \cite{Perelman}. 

In this paper, we consider similar questions in higher dimensions. We will assume throughout that $(M,g)$ is a steady gradient Ricci soliton of dimension $n \geq 4$ with positive sectional curvature. We may write $\text{\rm Ric} = D^2 f$ for some real-valued function $f$. As usual, we put $X = \nabla f$, and denote by $\Phi_t$ the flow generated by the vector field $-X$. 

\begin{definition} 
We say that $(M,g)$ is asymptotically cylindrical if the following holds: \\
(i) The scalar curvature satisfies $\frac{\Lambda_1}{d(p_0,p)} \leq R \leq \frac{\Lambda_2}{d(p_0,p)}$ at infinity, where $\Lambda_1$ and $\Lambda_2$ are positive constants. \\ 
(ii) Let $p_m$ be an arbitrary sequence of marked points going to infinity. Consider the rescaled metrics 
\[\hat{g}^{(m)}(t) = r_m^{-1} \, \Phi_{r_m t}^*(g),\] 
where $r_m \, R(p_m) = \frac{n-1}{2} + o(1)$. As $m \to \infty$, the flows $(M,\hat{g}^{(m)}(t),p_m)$ converge in the Cheeger-Gromov sense to a family of shrinking cylinders $(S^{n-1} \times \mathbb{R},\overline{g}(t))$, $t \in (0,1)$. The metric $\overline{g}(t)$ is given by 
\begin{equation} 
\label{shrinking.cylinders}
\overline{g}(t) = (n-2)(2-2t) \, g_{S^{n-1}} + dz \otimes dz, 
\end{equation}
where $g_{S^{n-1}}$ denotes the standard metric on $S^{n-1}$ with constant sectional curvature $1$. 
\end{definition}
We now state the main result of this paper. This result is motivated in part by the work of L.~Simon and B.~Solomon \cite{Simon-Solomon}, which deals with uniqueness questions for minimal surfaces with prescribed tangent cones at infinity. 

\begin{theorem} 
\label{main.theorem}
Let $(M,g)$ be a steady gradient Ricci soliton of dimension $n \geq 4$ which has positive sectional curvature and is asymptotically cylindrical. Then $(M,g)$ is rotationally symmetric. In particular, $(M,g)$ is isometric to the $n$-dimensional Bryant soliton up to scaling.
\end{theorem}

In dimension $3$, it follows from work of Perelman \cite{Perelman} that any complete steady gradient Ricci soliton which is non-flat and $\kappa$-noncollapsed is asymptotically cylindrical. Thus, Theorem \ref{main.theorem} can be viewed as a higher dimensional version of Theorem \ref{dim.3}. 

Theorem \ref{main.theorem} has an interesting implication in dimension $4$. A four-dimensional manifold $(M,g)$ has positive isotropic curvature if and only if $a_1+a_2 > 0$ and $c_1+c_2 > 0$, where $a_1,a_2,c_1,c_2$ are defined as in \cite{Hamilton1}. The notion of isotropic curvature was first introduced by Micallef and Moore \cite{Micallef-Moore} in their work on the index of minimal two-spheres. It also plays a central role in the convergence theory for the Ricci flow in higher dimensions (see e.g. \cite{Brendle1}, \cite{Brendle-book}).

\begin{theorem} 
\label{dim.4}
Let $(M,g)$ be a four-dimensional steady gradient Ricci soliton which is non-flat; is $\kappa$-noncollapsed; and satisfies the pointwise pinching condition 
\[0 \leq \max \{a_3,b_3,c_3\} \leq \Lambda \, \min \{a_1+a_2,c_1+c_2\},\] 
where $a_1,a_2,a_3,c_1,c_2,c_3,b_3$ are defined as in Hamilton's paper \cite{Hamilton1} and $\Lambda \geq 1$ is a constant. Then $(M,g)$ is rotationally symmetric.
\end{theorem}

We note that various authors have obtained uniqueness results for Ricci solitons in higher dimensions; see e.g. \cite{Cao-survey}, \cite{Cao-Chen}, \cite{Cao-Catino-Chen-Mantegazza-Mazzieri}, and \cite{Chen-Wang}. Moreover, T.~Ivey \cite{Ivey} has constructed examples of Ricci solitons which are not rotationally symmetric.

In order to prove Theorem \ref{main.theorem}, we will adapt the arguments in \cite{Brendle2}. While many arguments in \cite{Brendle2} directly generalize to higher dimensions, there are several crucial differences. In particular, the proof of the roundness estimate in Section \ref{roundness.estimate} is very different than in the three-dimensional case. Moreover, the proof in \cite{Brendle2} uses an estimate of Anderson and Chow \cite{Anderson-Chow} for the linearized Ricci flow system. This estimate uses special properties of the curvature tensor in dimension $3$, so we require a different argument to handle the higher dimensional case. This will be discussed in Section \ref{analysis.of.lich.eq}.

Finally, to deduce Theorem \ref{dim.4} from Theorem \ref{main.theorem}, we show that a steady gradient Ricci soliton $(M,g)$ which satisfies the assumptions of Theorem \ref{dim.4} must have positive curvature operator (cf. Corollary \ref{pco} below). The proof of this fact uses the pinching estimates of Hamilton (see \cite{Hamilton1}, \cite{Hamilton2}). Using results from \cite{Chen-Zhu}, we conclude that $(M,g)$ is asymptotically cylindrical. Theorem \ref{main.theorem} then implies that $(M,g)$ is rotationally symmetric.

\section{The roundness estimate}

\label{roundness.estimate}

By scaling, we may assume that $R + |\nabla f|^2 = 1$. Since $R \to 0$ at infinity, we can find a point $p_0$ where the scalar curvature attains its maximum. Since $(M,g)$ has positive sectional curvature, the Hessian of $f$ is strictly positive definite at each point in $M$. The identity $\nabla R(p_0) = 0$ implies $\nabla f(p_0) = 0$. Since $f$ is strictly convex, we conclude that $\liminf_{p \to \infty} \frac{f(p)}{d(p_0,p)} > 0$. On the other hand, since $|\nabla f|^2 \leq 1$, we have $\limsup_{p \to \infty} \frac{f(p)}{d(p_0,p)} < \infty$.

Using the fact that $(M,g)$ is asymptotically cylindrical, we obtain the following result:

\begin{proposition}
\label{asympt}
We have $f \, R = \frac{n-1}{2} + o(1)$ and $f \, \text{\rm Ric} \leq \big ( \frac{1}{2} + o(1) \big ) \, g$. Moreover, we have $f^2 \, \text{\rm Ric} \geq c \, g$ for some positive constant $c$.
\end{proposition}

\textbf{Proof.} 
Since $(M,g)$ is asymptotically cylindrical, we have $\Delta R = o(r^{-2})$ and $|\text{\rm Ric}|^2 = \frac{1}{n-1} \, R^2 + o(r^{-2})$. This implies 
\[-\langle X,\nabla R \rangle = \Delta R + 2 \, |\text{\rm Ric}|^2 = \frac{2}{n-1} \, R^2 + o(r^{-2}),\] 
hence 
\[\Big \langle X,\nabla \Big ( \frac{1}{R} - \frac{2}{n-1} \, f \Big ) \Big \rangle = o(1).\] 
Integrating this inequality along the integral curves of $X$ gives 
\[\frac{1}{R} - \frac{2}{n-1} \, f = o(r),\] 
hence 
\[f \, R = \frac{n-1}{2} + o(1).\] 
Moreover, we have $\text{\rm Ric} \leq \big ( \frac{1}{n-1} + o(1) \big ) \, R \, g$ since $(M,g)$ is asymptotically cylindrical. Therefore, $f \, \text{\rm Ric} \leq \big ( \frac{1}{2} + o(1) \big ) \, g$.

In order to verify the third statement, we choose an orthonormal frame $\{e_1,\hdots,e_n\}$ such that $e_n = \frac{X}{|X|}$. Since $(M,g)$ is asymptotically cylindrical, we have 
\[\text{\rm Ric}(e_i,e_j) = \frac{1}{n-1} \, R \, \delta_{ij} + o(r^{-1})\] 
for $i,j \in \{1,\hdots,n-1\}$ and 
\[2 \, \text{\rm Ric}(e_i,X) = -\langle e_i,\nabla R \rangle = o(r^{-\frac{3}{2}}).\] 
Moreover, we have 
\[2 \, \text{\rm Ric}(X,X) = -\langle X,\nabla R \rangle = \Delta R + 2 \, |\text{\rm Ric}|^2 = \frac{2}{n-1} \, R^2 + o(r^{-2}).\] 
Putting these facts together, we conclude that $\text{\rm Ric} \geq c \, R^2 \, g$ for some positive constant $c$. From this, the assertion follows. \\

In the remainder of this section, we prove a roundness estimate. We begin with a lemma: 

\begin{lemma}
\label{prelim.est}
We have $R_{ijkl} \, \partial^l f = O(r^{-\frac{3}{2}})$.
\end{lemma}

\textbf{Proof.} 
Using Shi's estimate, we obtain 
\[R_{ijkl} \, \partial^l f = D_i \text{\rm Ric}_{jk} - D_j \text{\rm Ric}_{ik} = O(r^{-\frac{3}{2}}).\] 
This proves the assertion. \\

We next define 
\[T = (n-1) \, \text{\rm Ric} - R \, g + R \, df \otimes df.\] 
Note that 
\begin{align*} 
&\text{\rm tr}(T) = -R^2 = O(r^{-2}), \\ 
&T(\nabla f,\cdot) = (n-1) \, \text{\rm Ric}(\nabla f,\cdot) - R^2 \, \nabla f = O(r^{-\frac{3}{2}}), \\ 
&T(\nabla f,\nabla f) = (n-1) \, \text{\rm Ric}(\nabla f,\nabla f) - R^2 \, |\nabla f|^2 = O(r^{-2}). 
\end{align*}

\begin{proposition} 
\label{roundness.1}
We have $|T| \leq O(r^{-\frac{3}{2}})$.
\end{proposition} 

\textbf{Proof.} 
The Ricci tensor of $(M,g)$ satisfies the equation  
\[\Delta \text{\rm Ric}_{ik} + D_X \text{\rm Ric}_{ik} = -2 \sum_{j,l=1}^n R_{ijkl} \, \text{\rm Ric}^{jl}.\] 
Moreover, using the identity $\Delta X + D_X X = 0$, we obtain  
\begin{align*} 
&\Delta (R \, g_{ik} - R \, \partial_i f \, \partial_k f) + D_X (R \, g_{ik} - R \, \partial_i f \, \partial_k f) \\ 
&= (\Delta R + \langle X,\nabla R \rangle) \, (g_{ik} - \partial_i f \, \partial_k f) + O(r^{-\frac{5}{2}}) \\ 
&= -2 \, |\text{\rm Ric}|^2 \, (g_{ik} - \partial_i f \, \partial_k f) + O(r^{-\frac{5}{2}}). 
\end{align*} 
Using Lemma \ref{prelim.est}, we conclude that 
\begin{align*} 
\Delta T_{ik} + D_X T_{ik} 
&= -2 \sum_{j,l=1}^{n-1} R_{ijkl} \, T^{jl} - 2 \, R \, \text{\rm Ric}_{ik} \\ 
&+ 2 \, |\text{\rm Ric}|^2 \, (g_{ik} - \partial_i f \, \partial_k f) + O(r^{-\frac{5}{2}}), 
\end{align*}
hence 
\begin{align*} 
&\Delta(|T|^2) + \langle X,\nabla(|T|^2) \rangle \\ 
&= 2 \, |DT|^2 - 4 \sum_{j,l=1}^{n-1} R_{ijkl} \, T^{ik} \, T^{jl} - 4 \, R \sum_{i,k=1}^n \text{\rm Ric}_{ik} \, T^{ik} \\ 
&+ 4 \, |\text{\rm Ric}|^2 \sum_{i,k=1}^n (g_{ik} - \partial_i f \, \partial_k f) \, T^{ik} + O(r^{-\frac{5}{2}}) \, |T| \\ 
&= 2 \, |DT|^2 - 4 \sum_{j,l=1}^{n-1} R_{ijkl} \, T^{ik} \, T^{jl} - \frac{4}{n-1} \, R \, |T|^2 \\ 
&+ 4 \, \Big ( |\text{\rm Ric}|^2 - \frac{1}{n-1} \, R^2 \Big ) \, \sum_{i,k=1}^n (g_{ik} - \partial_i f \, \partial_k f) \, T^{ik} + O(r^{-\frac{5}{2}}) \, |T|. 
\end{align*} 
Since $\sum_{i,k=1}^n (g_{ik} - \partial_i f \, \partial_k f) \, T^{ik} = O(r^{-2})$, we obtain 
\begin{align*} 
&\Delta(|T|^2) + \langle X,\nabla(|T|^2) \rangle \\ 
&\geq - 4 \sum_{j,l=1}^{n-1} R_{ijkl} \, T^{ik} \, T^{jl} - \frac{4}{n-1} \, R \, |T|^2 - O(r^{-\frac{5}{2}}) \, |T| - O(r^{-4}). 
\end{align*} 
Moreover, since $(M,g)$ is asymptotically cylindrical, we have 
\begin{align*} 
R_{ijkl} 
&= \frac{1}{(n-1)(n-2)} \, R \, (g_{ik} - \partial_i f \, \partial_k f) \, (g_{jl} - \partial_j f \, \partial_l f) \\ 
&- \frac{1}{(n-1)(n-2)} \, R \, (g_{il} - \partial_i f \, \partial_l f) \, (g_{jk} - \partial_j f \, \partial_k f) \\ 
&+ o(r^{-1}) 
\end{align*}
near infinity. This implies 
\[\sum_{j,l=1}^{n-1} R_{ijkl} \, T^{ik} \, T^{jl} = -\frac{1}{(n-1)(n-2)} \, R \, |T|^2 + O(r^{-\frac{5}{2}}) \, |T| + o(r^{-1}) \, |T|^2,\] 
hence 
\begin{align*} 
&\Delta(|T|^2) + \langle X,\nabla(|T|^2) \rangle \\ 
&\geq -\frac{4(n-3)}{(n-1)(n-2)} \, R \, |T|^2 - o(r^{-1}) \, |T|^2 - O(r^{-\frac{5}{2}}) \, |T| - O(r^{-4}). 
\end{align*}
We next observe that $|D_X \text{\rm Ric}| \leq O(r^{-2})$ and $|D_{X,X}^2 \text{\rm Ric}| \leq O(r^{-\frac{5}{2}})$. This implies $|D_X T| \leq O(r^{-2})$ and $|D_{X,X}^2 T| \leq O(r^{-\frac{5}{2}})$. From this, we deduce that 
\begin{align*} 
&\Delta_\Sigma(|T|^2) + \langle X,\nabla(|T|^2) \rangle \\ 
&\geq -\frac{2(n-3)}{n-2} \, f^{-1} \, |T|^2 - o(r^{-1}) \, |T|^2 - O(r^{-\frac{5}{2}}) \, |T| - O(r^{-4}), 
\end{align*}
where $\Delta_\Sigma$ denotes the Laplacian on the level surfaces of $f$. Thus, we conclude that 
\begin{align*} 
&\Delta_\Sigma(f^2 \, |T|^2) + \langle X,\nabla(f^2 \, |T|^2) \rangle \\ 
&\geq \frac{2}{n-2} \, f \, |T|^2 - o(r) \, |T|^2 - O(r^{-\frac{1}{2}}) \, |T| - O(r^{-2}) \geq -O(r^{-2}) 
\end{align*}
outside some compact set. Since $f^2 \, |T|^2 \to 0$ at infinity, the parabolic maximum principle implies that $f^2 \, |T|^2 \leq O(r^{-1})$. This completes the proof. \\

In the following, we fix $\varepsilon$ sufficiently small; for example, $\varepsilon = \frac{1}{1000n}$ will work. By Proposition \ref{roundness.1}, we have $|T| \leq O(r^{\frac{1}{2(n-2)} - \frac{3}{2} - 32\varepsilon})$. Moreover, it follows from Shi's estimates that $|D^m T| \leq O(r^{-\frac{m+2}{2}})$ for each $m$. Using standard interpolation inequalities, we obtain $|DT| \leq O(r^{\frac{1}{2(n-2)}-2-16\varepsilon})$. Using the identity 
\begin{align*} 
D^k T_{ik} 
&= \frac{n-3}{2} \, \partial_i R + \langle \nabla f,\nabla R \rangle \, \partial_i f + R^2 \, \partial_i f + R \, \text{\rm Ric}_i^k \, \partial_k f \\ 
&= \frac{n-3}{2} \, \partial_i R + O(r^{-2}), 
\end{align*}
we conclude that $|\nabla R| \leq O(r^{\frac{1}{2(n-2)}-2-16\varepsilon})$. This implies 
\[|D\text{\rm Ric}| \leq C \, |DT| + C \, |\nabla R| + C \, R \, |D^2 f| \leq O(r^{\frac{1}{2(n-2)}-2-16\varepsilon}).\] 
Using standard interpolation inequalities, we obtain 
\[|D^2 \text{\rm Ric}| \leq O(r^{\frac{1}{2(n-2)}-\frac{5}{2}-8\varepsilon}).\] 

\begin{proposition}
\label{scalar.curvature}
We have $f \, R = \frac{n-1}{2} + O(r^{\frac{1}{2(n-2)}-\frac{1}{2}-8\varepsilon})$.
\end{proposition} 

\textbf{Proof.} 
Using the inequality $|T| \leq O(r^{-\frac{3}{2}})$, we obtain 
\[|\text{\rm Ric}| = \frac{1}{n-1} \ R \, |g - df \otimes df| + O(r^{-\frac{3}{2}}) = \frac{1}{\sqrt{n-1}} \, R + O(r^{-\frac{3}{2}}),\] 
hence 
\[|\text{\rm Ric}|^2 = \frac{1}{n-1} \, R^2 + O(r^{-\frac{5}{2}}).\] 
This implies 
\[-\langle X,\nabla R \rangle = \Delta R + 2 \, |\text{\rm Ric}|^2 = \frac{2}{n-1} \, R^2 + O(r^{\frac{1}{2(n-2)}-\frac{5}{2}-8\varepsilon}),\] 
hence 
\[\Big \langle X,\nabla \Big ( \frac{1}{R} - \frac{2}{n-1} \, f \Big ) \Big \rangle = O(r^{\frac{1}{2(n-2)}-\frac{1}{2}-8\varepsilon}).\] 
Integrating this identity along the integral curves of $X$, we obtain 
\[\frac{1}{R} - \frac{2}{n-1} \, f = O(r^{\frac{1}{2(n-2)}+\frac{1}{2}-8\varepsilon}).\] From this, the assertion follows. \\

\begin{proposition} 
\label{roundness.2}
We have 
\begin{align*} 
f \, R_{ijkl} 
&= \frac{1}{2(n-2)} \, (g_{ik} - \partial_i f \, \partial_k f) \, (g_{ik} - \partial_i f \, \partial_k f) \\ 
&- \frac{1}{2(n-2)} \, (g_{il} - \partial_i f \, \partial_l f) \, (g_{jk} - \partial_j f \, \partial_k f) \\ 
&+ O(r^{\frac{1}{2(n-2)}-\frac{1}{2}-8\varepsilon}). 
\end{align*}
\end{proposition}

\textbf{Proof.} 
It follows from Proposition 2.10 in \cite{Brendle-book} that 
\begin{align*} 
-D_X R_{ijkl} 
&= D_{i,k}^2 \text{\rm Ric}_{jl} - D_{i,l}^2 \text{\rm Ric}_{jk} - D_{j,k}^2 \text{\rm Ric}_{il} + D_{j,l}^2 \text{\rm Ric}_{ik} \\ 
&+ \sum_{m=1}^n \text{\rm Ric}_i^m \, R_{mjkl} + \sum_{m=1}^n \text{\rm Ric}_j^m \, R_{imkl}. 
\end{align*}
Using Lemma \ref{prelim.est} and Proposition \ref{roundness.1}, we obtain 
\begin{align*} 
\sum_{m=1}^n \text{\rm Ric}_i^m \, R_{mjkl} 
&= \frac{1}{n-1} \, R \sum_{m=1}^n (\delta_i^m - \partial_i f \, \partial^m f) \, R_{mjkl} + O(r^{-\frac{5}{2}}) \\ 
&= \frac{1}{n-1} \, R \, R_{ijkl} + O(r^{-\frac{5}{2}}). 
\end{align*}
Thus, we conclude that 
\begin{align*} 
-D_X R_{ijkl} 
&= \frac{2}{n-1} \, R \, R_{ijkl} + O(r^{\frac{1}{2(n-2)}-\frac{5}{2}-8\varepsilon}) \\ 
&= f^{-1} \, R_{ijkl} + O(r^{\frac{1}{2(n-2)}-\frac{5}{2}-8\varepsilon}), 
\end{align*}
hence 
\[|D_X(f \, R_{ijkl})| \leq O(r^{\frac{1}{2(n-2)}-\frac{3}{2}-8\varepsilon}).\] 
On the other hand, the tensor 
\begin{align*} 
S_{ijkl} 
&= \frac{1}{2(n-2)} \, (g_{ik} - \partial_i f \, \partial_k f) \, (g_{jl} - \partial_j f \, \partial_l f) \\ 
&- \frac{1}{2(n-2)} \, (g_{il} - \partial_i f \, \partial_l f) \, (g_{jk} - \partial_j f \, \partial_k f) 
\end{align*} 
satisfies 
\[|D_X S_{ijkl}| \leq O(r^{-\frac{3}{2}}).\] 
Putting these facts together, we obtain 
\[|D_X(f \, R_{ijkl} - S_{ijkl})| \leq O(r^{\frac{1}{2(n-2)}-\frac{3}{2}-8\varepsilon}).\] 
Moreover, we have $|f \, R_{ijkl} - S_{ijkl}| \to 0$ at infinity. Integrating the preceding inequality along integral curves of $X$ gives 
\[|f \, R_{ijkl} - S_{ijkl}| \leq O(r^{\frac{1}{2(n-2)}-\frac{1}{2}-8\varepsilon}),\] 
as claimed. \\

We next construct a collection of approximate Killing vector fields: 

\begin{proposition}
\label{rotation.vector.fields}
We can find a collection of vector fields $U_a$, $a \in \{1,\hdots,\frac{n(n-1)}{2}\}$, on $(M,g)$ such that $|\mathscr{L}_{U_a}(g)| \leq O(r^{\frac{1}{2(n-2)}-\frac{1}{2}-2\varepsilon})$ and $|\Delta U_a + D_X U_a| \leq O(r^{\frac{1}{2(n-2)}-1-2\varepsilon})$. Moreover, we have 
\[\sum_{a=1}^{\frac{n(n-1)}{2}} U_a \otimes U_a = r \, \Big ( \sum_{i=1}^{n-1} e_i \otimes e_i + O(r^{\frac{1}{2(n-2)}-\frac{1}{2}-2\varepsilon}) \Big ),\] 
where $\{e_1,\hdots,e_{n-1}\}$ is a local orthonormal frame on the level set $\{f=r\}$.
\end{proposition}

The proof of Proposition \ref{rotation.vector.fields} is analogous to the arguments in \cite{Brendle2}, Section 3. We omit the details.

\section{An elliptic PDE for vector fields}

\label{pde.for.vector.fields}

Let us fix a smooth vector field $Q$ on $M$ with the property that $|Q| \leq O(r^{\frac{1}{2(n-2)}-1-2\varepsilon})$. We will show that there exists a vector field $V$ on $M$ such that $\Delta V + D_X V = Q$ and $|V| \leq O(r^{\frac{1}{2(n-2)}-\varepsilon})$. 

\begin{lemma} 
\label{spherical.harmonics.V}
Consider the shrinking cylinders $(S^{n-1} \times \mathbb{R},\overline{g}(t))$, $t \in (0,1)$, where $\overline{g}(t)$ is given by (\ref{shrinking.cylinders}). Let $\overline{V}(t)$, $t \in (0,1)$, be a one-parameter family of vector fields which satisfy the parabolic equation 
\begin{equation} 
\label{evolution.of.vector.field}
\frac{\partial}{\partial t} \overline{V}(t) = \Delta_{\overline{g}(t)} \overline{V}(t) + \text{\rm Ric}_{\overline{g}(t)}(\overline{V}(t)). 
\end{equation}
Moreover, suppose that $\overline{V}(t)$ is invariant under translations along the axis of the cylinder, and 
\begin{equation}
\label{boundedness.assumption.V}
|\overline{V}(t)|_{\overline{g}(t)} \leq 1 
\end{equation} 
for all $t \in (0,\frac{1}{2}]$. Then 
\[\inf_{\lambda \in \mathbb{R}} \sup_{S^{n-1} \times \mathbb{R}} \Big | \overline{V}(t) - \lambda \, \frac{\partial}{\partial z} \Big |_{\overline{g}(t)} \leq L \, (1-t)^{\frac{1}{2(n-2)}}\] 
for all $t \in [\frac{1}{2},1)$, where $L$ is a positive constant.
\end{lemma}

\textbf{Proof.} 
Since $\overline{V}(t)$ is invariant under translations along the axis of the cylinder, we may write 
\[\overline{V}(t) = \xi(t) + \eta(t) \, \frac{\partial}{\partial z}\] 
for $t \in (0,1)$, where $\xi(t)$ is a vector field on $S^{n-1}$ and $\eta(t)$ is a real-valued function on $S^{n-1}$. The parabolic equation (\ref{evolution.of.vector.field}) implies the following system of equations for $\xi(t)$ and $\eta(t)$: 
\begin{align} 
&\frac{\partial}{\partial t} \xi(t) = \frac{1}{(n-2)(2-2t)} \, (\Delta_{S^{n-1}} \xi(t) + (n-2) \, \xi(t)), \label{vec.1} \\ 
&\frac{\partial}{\partial t} \eta(t) = \frac{1}{(n-2)(2-2t)} \, \Delta_{S^{n-1}} \eta(t). \label{vec.2} 
\end{align}
Furthermore, the estimate (\ref{boundedness.assumption.V}) gives 
\begin{align} 
&\sup_{S^{n-1}} |\xi(t)|_{g_{S^{n-1}}} \leq L_1, \label{est.1} \\ 
&\sup_{S^{n-1}} |\eta(t)| \leq L_1 \label{est.2}
\end{align} 
for each $t \in (0,\frac{1}{2}]$, where $L_1$ is a positive constant.

Let us consider the operator $\xi \mapsto -\Delta_{S^{n-1}} \xi - (n-2) \, \xi$, acting on vector fields on $S^{n-1}$. By Proposition \ref{laplacian.on.one.forms}, the first eigenvalue of this operator is at least $-(n-3)$. Using (\ref{vec.1}) and (\ref{est.1}), we obtain 
\begin{equation} 
\label{est.3} 
\sup_{S^{n-1}} |\xi(t)|_{g_{S^{n-1}}} \leq L_2 \, (1-t)^{-\frac{n-3}{2(n-2)}}
\end{equation} 
for all $t \in [\frac{1}{2},1)$, where $L_2$ is a positive constant. Similarly, it follows from (\ref{vec.2}) and (\ref{est.2}) that 
\begin{equation} 
\label{est.4}
\inf_{\lambda \in \mathbb{R}} \sup_{S^{n-1}} |\eta(t) - \lambda| \leq L_3 \, (1-t)^{\frac{n-1}{2(n-2)}} 
\end{equation} 
for each $t \in [\frac{1}{2},1)$, where $L_3$ is a positive constant. Combining (\ref{est.3}) and (\ref{est.4}), the assertion follows. \\

\begin{lemma}[cf. \cite{Brendle2}, Lemma 5.2]
\label{max.principle.V}
Let $V$ be a smooth vector field satisfying $\Delta V + D_X V = Q$ in the region $\{f \leq \rho\}$. Then 
\[\sup_{\{f \leq \rho\}} |V| \leq \sup_{\{f=\rho\}} |V| + B \, \rho^{\frac{1}{2(n-2)}-2\varepsilon}\] 
for some uniform constant $B \geq 1$.
\end{lemma} 

The proof of Lemma \ref{max.principle.V} is similar to the proof of Lemma 5.2 in \cite{Brendle2}; we omit the details.

As in \cite{Brendle2}, we choose a sequence of real numbers $\rho_m \to \infty$. For each $m$, we can find a vector field $V^{(m)}$ such that $\Delta V^{(m)} + D_X V^{(m)} = Q$ in the region $\{f \leq \rho_m\}$ and $V^{(m)} = 0$ on the boundary $\{f=\rho_m\}$. We now define  
\[A^{(m)}(r) = \inf_{\lambda \in \mathbb{R}} \sup_{\{f=r\}} |V^{(m)} - \lambda \, X|\] 
for $r \leq \rho_m$.

\begin{lemma}
\label{iteration}
Let us fix a real number $\tau \in (0,\frac{1}{2})$ so that $\tau^{-\varepsilon} > 2L$, where $L$ is the constant in Lemma \ref{spherical.harmonics.V}. Then we can find a real number $\rho_0$ and a positive integer $m_0$ such that 
\[2 \, \tau^{-\frac{1}{2(n-2)}+\varepsilon} \, A^{(m)}(\tau r) \leq A^{(m)}(r) + r^{\frac{1}{2(n-2)}-\varepsilon}\] 
for all $r \in [\rho_0,\rho_m]$ and all $m \geq m_0$.
\end{lemma}

\textbf{Proof.} 
We argue by contradiction. Suppose that the assertion is false. After passing to a subsequence, there exists a sequence of real numbers $r_m \leq \rho_m$ such that $r_m \to \infty$ and 
\[A^{(m)}(r_m) + r_m^{\frac{1}{2(n-2)}-\varepsilon} \leq 2 \, \tau^{-\frac{1}{2(n-2)}+\varepsilon} \, A^{(m)}(\tau r_m)\] 
for all $m$. For each $m$, there exists a real number $\lambda_m$ such that 
\[\sup_{\{f=r_m\}} |V^{(m)} - \lambda_m \, X| = A^{(m)}(r_m).\] 
Applying Lemma \ref{max.principle.V} to the vector field $V^{(m)} - \lambda_m \, X$ gives 
\begin{align*} 
\sup_{\{f \leq r_m\}} |V^{(m)} - \lambda_m \, X| 
&\leq \sup_{\{f=r_m\}} |V^{(m)} - \lambda_m \, X| + B \, r_m^{\frac{1}{2(n-2)}-2\varepsilon} \\ 
&\leq A^{(m)}(r_m) + r_m^{\frac{1}{2(n-2)}-\varepsilon} 
\end{align*} 
if $m$ is sufficiently large. We next consider the vector field 
\[\tilde{V}^{(m)} = \frac{1}{A^{(m)}(r_m) + r_m^{\frac{1}{2(n-2)}-\varepsilon}} \, (V^{(m)} - \lambda_m \, X).\] 
The vector field $\tilde{V}^{(m)}$ satisfies 
\begin{equation} 
\label{bound.for.V}
\sup_{\{f \leq r_m\}} |\tilde{V}^{(m)}| \leq 1. 
\end{equation} 
Let 
\[\hat{g}^{(m)}(t) = r_m^{-1} \, \Phi_{r_m t}^*(g)\] 
and 
\[\hat{V}^{(m)}(t) = r_m^{\frac{1}{2}} \, \Phi_{r_m t}^*(\tilde{V}^{(m)}).\] 
Note that the metrics $\hat{g}^{(m)}(t)$ evolve by the Ricci flow. Moreover, the vector fields $\hat{V}^{(m)}(t)$ satisfy the parabolic equation 
\[\frac{\partial}{\partial t} \hat{V}^{(m)}(t) = \Delta_{\hat{g}^{(m)}(t)} \hat{V}^{(m)}(t) + \text{\rm Ric}_{\hat{g}^{(m)}(t)}(\hat{V}^{(m)}(t)) - \hat{Q}^{(m)}(t),\] 
where 
\[\hat{Q}^{(m)}(t) = \frac{r_m^{\frac{3}{2}}}{A^{(m)}(r_m) + r_m^{\frac{1}{2(n-2)}-\varepsilon}} \, \Phi_{r_mt}^*(Q).\] 
Using (\ref{bound.for.V}), we obtain 
\[\limsup_{m \to \infty} \sup_{t \in [\delta,1-\delta]} \sup_{\{r_m - \delta^{-1} \, \sqrt{r_m} \leq f \leq r_m + \delta^{-1} \, \sqrt{r_m}\}} |\hat{V}^{(m)}(t)|_{\hat{g}^{(m)}(t)} < \infty\] 
for any given $\delta \in (0,\frac{1}{2})$. Moreover, the estimate $|Q| \leq O(r^{\frac{1}{2(n-2)}-1-2\varepsilon})$ implies that 
\[\limsup_{m \to \infty} \sup_{t \in [\delta,1-\delta]} \sup_{\{r_m - \delta^{-1} \, \sqrt{r_m} \leq f \leq r_m + \delta^{-1} \, \sqrt{r_m}\}} |\hat{Q}^{(m)}(t)|_{\hat{g}^{(m)}(t)} = 0\] 
for any given $\delta \in (0,\frac{1}{2})$. 

We now pass to the limit as $m \to \infty$. To that end, we choose a sequence of marked points $p_m \in M$ such that $f(p_m) = r_m$. The manifolds $(M,\hat{g}^{(m)}(t),p_m)$ converge in the Cheeger-Gromov sense to a one-parameter family of shrinking cylinders $(S^{n-1} \times \mathbb{R},\overline{g}(t))$, $t \in (0,1)$, where $\overline{g}(t)$ is given by (\ref{shrinking.cylinders}). Furthermore, the rescaled vector fields $r_m^{\frac{1}{2}} \, X$ converge to the axial vector field $\frac{\partial}{\partial z}$ on $S^{n-1} \times \mathbb{R}$. Finally, the sequence $\hat{V}^{(m)}(t)$ converges in $C_{loc}^0$ to a one-parameter family of vector fields $\overline{V}(t)$, $t \in (0,1)$, which satisfy the parabolic equation 
\[\frac{\partial}{\partial t} \overline{V}(t) = \Delta_{\overline{g}(t)} \overline{V}(t) + \text{\rm Ric}_{\overline{g}(t)}(\overline{V}(t)).\] 
As in \cite{Brendle2}, we can show that $\overline{V}(t)$ is invariant under translations along the axis of the cylinder. Moreover, the estimate (\ref{bound.for.V}) implies that 
\[|\overline{V}(t)|_{\overline{g}(t)} \leq 1\] 
for all $t \in (0,\frac{1}{2}]$. Hence, it follows from Lemma \ref{spherical.harmonics.V} that 
\begin{equation} 
\label{symmetry.V}
\inf_{\lambda \in \mathbb{R}} \sup_{S^{n-1} \times \mathbb{R}} \Big | \overline{V}(t) - \lambda \, \frac{\partial}{\partial z} \Big |_{\overline{g}(t)} \leq L \, (1-t)^{\frac{1}{2(n-2)}} 
\end{equation}
for all $t \in (0,\frac{1}{2}]$. Finally, we have 
\begin{align*} 
&\inf_{\lambda \in \mathbb{R}} \sup_{\Phi_{r_m(\tau-1)}(\{f=\tau r_m\})} \Big | \hat{V}^{(m)}(1-\tau) - \lambda \, r_m^{\frac{1}{2}} \, X \Big |_{\hat{g}^{(m)}(1-\tau)} \\ 
&= \inf_{\lambda \in \mathbb{R}} \sup_{\{f=\tau r_m\}} |\tilde{V}^{(m)} - \lambda \, X|_g \\ 
&= \frac{1}{A^{(m)}(r_m)+ r_m^{\frac{1}{2(n-2)}-\varepsilon}} \, \inf_{\lambda \in \mathbb{R}} \sup_{\{f=\tau r_m\}} |V^{(m)} - \lambda \, X|_g \\ 
&= \frac{A^{(m)}(\tau r_m)}{A^{(m)}(r_m) + r_m^{\frac{1}{2(n-2)}-\varepsilon}} \\ 
&\geq \frac{1}{2} \, \tau^{\frac{1}{2(n-2)}-\varepsilon}. 
\end{align*} 
If we send $m \to \infty$, we obtain 
\begin{equation} 
\label{asymmetry.V}
\inf_{\lambda \in \mathbb{R}} \sup_{S^{n-1} \times \mathbb{R}} \Big | \overline{V}(1-\tau) - \lambda \, \frac{\partial}{\partial z} \Big |_{\overline{g}(1-\tau)} \geq \frac{1}{2} \, \tau^{\frac{1}{2(n-2)}-\varepsilon}. 
\end{equation} 
Since $\tau^{-\varepsilon} > 2L$, the inequality (\ref{asymmetry.V}) is in contradiction with (\ref{symmetry.V}). This completes the proof of Lemma \ref{iteration}. \\

If we iterate the estimate in Lemma \ref{iteration}, we obtain 
\[\sup_m \sup_{\rho_0 \leq r \leq \rho_m} r^{-\frac{1}{2(n-2)}+\varepsilon} \, A^{(m)}(r) < \infty.\] From this, we deduce the following result:

\begin{proposition} 
\label{uniform.bounds}
There exists a sequence of real numbers $\lambda_m$ such that 
\[\sup_m \sup_{\{f \leq \rho_m\}} f^{-\frac{1}{2(n-2)}+\varepsilon} \, |V^{(m)} - \lambda_m \, X| < \infty.\] 
\end{proposition}

The proof of Proposition \ref{uniform.bounds} is analogous to the proof of Proposition 5.4 in \cite{Brendle2}. We omit the details. By taking the limit as $m \to \infty$ of the vector fields $V^{(m)} - \lambda_m \, X$, we obtain the following result:

\begin{theorem} 
\label{existence.of.V}
There exists a smooth vector field $V$ such that $\Delta V + D_X V = Q$ and $|V| \leq O(r^{\frac{1}{2(n-2)}-\varepsilon})$. Moreover, $|DV| \leq O(r^{\frac{1}{2(n-2)}-\frac{1}{2}-\varepsilon})$.
\end{theorem}

\section{Analysis of the Lichnerowicz equation}

\label{analysis.of.lich.eq}

Throughout this section, we will denote by $\Delta_L$ the Lichnerowicz Laplacian; that is, 
\[\Delta_L h_{ik} = \Delta h_{ik} + 2 \, R_{ijkl} \, h^{jl} - \text{\rm Ric}_i^l \, h_{kl} - \text{\rm Ric}_k^l \, h_{il}.\]

\begin{lemma} 
\label{spherical.harmonics.h}
Let us consider the shrinking cylinders $(S^{n-1} \times \mathbb{R},\overline{g}(t))$, $t \in (0,1)$, where $\overline{g}(t)$ is given by (\ref{shrinking.cylinders}). Let $\overline{h}(t)$, $t \in (0,1)$, be a one-parameter family of $(0,2)$-tensors which solve the parabolic equation 
\begin{equation} 
\label{parabolic.lichnerowicz.equation}
\frac{\partial}{\partial t} \overline{h}(t) = \Delta_{L,\overline{g}(t)} \overline{h}(t). 
\end{equation} 
Moreover, suppose that $\bar{h}(t)$ is invariant under translations along the axis of the cylinder, and 
\begin{equation} 
\label{boundedness.assumption.h}
|\overline{h}(t)|_{\overline{g}(t)} \leq (1-t)^{-2} 
\end{equation} 
for all $t \in (0,\frac{1}{2}]$. Then 
\[\inf_{\lambda \in \mathbb{R}} \sup_{S^{n-1} \times \mathbb{R}} \big | \overline{h}(t) - \lambda \, \text{\rm Ric}_{\overline{g}(t)} \big |_{\overline{g}(t)} \leq N \, (1-t)^{\frac{1}{2(n-2)}-\frac{1}{2}}\] 
for all $t \in [\frac{1}{2},1)$, where $N$ is a positive constant.
\end{lemma}

\textbf{Proof.} 
Since $\overline{h}(t)$ is invariant under translations along the axis of the cylinder, we may write 
\[\overline{h}(t) = \chi(t) + dz \otimes \sigma(t) + \sigma(t) \otimes dz + \beta(t) \, dz \otimes dz\] 
for $t \in (0,1)$, where $\chi(t)$ is a symmetric $(0,2)$ tensor on $S^{n-1}$, $\sigma(t)$ is a one-form on $S^{n-1}$, and $\beta(t)$ is a real-valued function on $S^{n-1}$. The parabolic Lichnerowicz equation (\ref{parabolic.lichnerowicz.equation}) implies the following system of equations for $\chi(t)$, $\sigma(t)$, and $\beta(t)$: 
\begin{align} 
&\frac{\partial}{\partial t} \chi(t) = \frac{1}{(n-2)(2-2t)} \, (\Delta_{S^{n-1}} \chi(t) - 2(n-1) \, \overset{\text{\rm o}}{\chi}(t)), \label{lich.1} \\ 
&\frac{\partial}{\partial t} \sigma(t) = \frac{1}{(n-2)(2-2t)} \, (\Delta_{S^{n-1}} \sigma(t) - (n-2) \, \sigma(t)), \label{lich.2} \\ 
&\frac{\partial}{\partial t} \beta(t) = \frac{1}{(n-2)(2-2t)} \, \Delta_{S^{n-1}} \beta(t). \label{lich.3}
\end{align}
Here, $\overset{\text{\rm o}}{\chi}(t)$ denotes the trace-free part of $\chi(t)$ with respect to the standard metric on $S^{n-1}$. Using the assumption (\ref{boundedness.assumption.h}), we obtain
\begin{align} 
&\sup_{S^{n-1}} |\chi(t)|_{g_{S^{n-1}}} \leq N_1, \label{estimate.1} \\ 
&\sup_{S^{n-1}} |\sigma(t)|_{g_{S^{n-1}}} \leq N_1, \label{estimate.2} \\ 
&\sup_{S^{n-1}} |\beta(t)| \leq N_1 \label{estimate.3}
\end{align} 
for each $t \in (0,\frac{1}{2}]$, where $N_1$ is a positive constant.

We next analyze the operator $\chi \mapsto -\Delta_{S^{n-1}} \chi + 2(n-1) \, \overset{\text{\rm o}}{\chi}$, acting on symmetric $(0,2)$-tensors on $S^{n-1}$. The first eigenvalue of this operator is equal to $0$, and the associated eigenspace is spanned by $g_{S^{n-1}}$. Moreover, the other eigenvalues of this operator are at least $n-1$ (cf. Proposition \ref{laplacian.on.tensors} below). Hence, it follows from (\ref{lich.1}) and (\ref{estimate.1}) that 
\begin{equation} 
\label{estimate.4}
\inf_{\lambda \in \mathbb{R}} \sup_{S^{n-1}} |\chi(t) - \lambda \, g_{S^{n-1}}|_{g_{S^{n-1}}} \leq N_2 \, (1-t)^{\frac{n-1}{2(n-2)}} 
\end{equation} 
for all $t \in [\frac{1}{2},1)$, where $N_2$ is a positive constant. We now consider the operator $\sigma \mapsto -\Delta_{S^{n-1}} \sigma + (n-2)\sigma$, acting on one-forms on $S^{n-1}$. By Proposition \ref{laplacian.on.one.forms}, the first eigenvalue of this operator is at least $n-1$. Using (\ref{lich.2}) and (\ref{estimate.2}), we deduce that 
\begin{equation}
\label{estimate.5}
\sup_{S^{n-1}} |\sigma(t)|_{g_{S^{n-1}}} \leq N_3 \, (1-t)^{\frac{n-1}{2(n-2)}} 
\end{equation}
for all $t \in [\frac{1}{2},1)$, where $N_3$ is a positive constant. Finally, using (\ref{lich.3}) and (\ref{estimate.3}), we obtain 
\begin{equation} 
\label{estimate.6}
\sup_{S^{n-1}} |\beta(t)| \leq N_4 
\end{equation} 
for all $t \in [\frac{1}{2},1)$, where $N_4$ is a positive constant. If we combine (\ref{estimate.4}), (\ref{estimate.5}), and (\ref{estimate.6}), the assertion follows. \\

We now study the equation $\Delta_L h + \mathscr{L}_X(h) = 0$ on $(M,g)$, where $\Delta_L$ denotes the Lichnerowicz Laplacian defined above.

\begin{lemma}
\label{prep}
Let $h$ be a solution of the Lichnerowicz-type equation 
\[\Delta_L h + \mathscr{L}_X(h) = 0\] 
on the region $\{f \leq \rho\}$. Then 
\[\sup_{\{f \leq \rho\}} |h| \leq C \, \rho^2 \, \sup_{\{f=\rho\}} |h|\] 
for some uniform constant $C$ which is independent of $\rho$.
\end{lemma}

\textbf{Proof.} 
It suffices to show that 
\begin{equation} 
\label{upper.bound}
h \leq C \, \rho^2 \, \Big ( \sup_{\{f=\rho\}} |h| \Big ) \, g 
\end{equation}
for some uniform constant $C$. Indeed, if (\ref{upper.bound}) holds, the assertion follows by applying (\ref{upper.bound}) to $h$ and $-h$.

We now describe the proof of (\ref{upper.bound}). By Proposition \ref{asympt}, we have $f^2 \, \text{\rm Ric} \geq c \, g$ for some positive constant $c$. Therefore, the tensor $\text{\rm Ric} - \frac{c}{2} \, \rho^{-2} \, g$ is positive definite in the region $\{f \leq \rho\}$. Let $\theta$ be the smallest real number with the property that $\theta \, (\text{\rm Ric} - \frac{c}{2} \, \rho^{-2} \, g) - h$ is positive semi-definite at each point in the region $\{f \leq \rho\}$. There exists a point $p_0 \in \{f \leq \rho\}$ and an orthonormal basis $\{e_1,\hdots,e_n\}$ of $T_{p_0} M$ such that 
\[\theta \, \text{\rm Ric}(e_1,e_1) - \frac{\theta \, c}{2} \, \rho^{-2} - h(e_1,e_1) = 0\] 
at the point $p_0$. We now distinguish two cases: 

\textit{Case 1:} Suppose that $p_0 \in \{f<\rho\}$. In this case, we have 
\[\theta \, (\Delta \text{\rm Ric})(e_1,e_1) - (\Delta h)(e_1,e_1) \geq 0\] 
and 
\[\theta \, (D_X \text{\rm Ric})(e_1,e_1) - (D_X h)(e_1,e_1) = 0\] 
at the point $p_0$. Using the identity $\Delta_L h + \mathscr{L}_X(h) = 0$, we obtain 
\begin{align*} 
0 &= (\Delta h)(e_1,e_1) + (D_X h)(e_1,e_1) + 2 \sum_{i,k=1}^n R(e_1,e_i,e_1,e_k) \, h(e_i,e_k) \\ 
&\leq \theta \, (\Delta \text{\rm Ric})(e_1,e_1) + \theta \, (D_X \text{\rm Ric})(e_1,e_1) + 2 \sum_{i,k=1}^n R(e_1,e_i,e_1,e_k) \, h(e_i,e_k) \\ 
&= -2 \sum_{i,k=1}^n R(e_1,e_i,e_1,e_k) \, (\theta \, \text{\rm Ric}(e_i,e_k) - h(e_i,e_k)) \\ 
&= -\theta \, c \, \rho^{-2} \, \text{\rm Ric}(e_1,e_1) \\ 
&- 2 \sum_{i,k=1}^n R(e_1,e_i,e_1,e_k) \, \Big ( \theta \, \text{\rm Ric}(e_i,e_k) - \frac{\theta \, c}{2} \, \rho^{-2} \, g(e_i,e_k) - h(e_i,e_k) \Big )
\end{align*} 
at the point $p_0$. Since $(M,g)$ has positive sectional curvature, we have 
\[\sum_{i,k=1}^n R(e_1,e_i,e_1,e_k) \, \Big ( \theta \, \text{\rm Ric}(e_i,e_k) - \frac{\theta \, c}{2} \, \rho^{-2} \, g(e_i,e_k) - h(e_i,e_k) \Big ) \geq 0.\] 
Consequently, $\theta \leq 0$. This implies $h \leq 0$ at each point in the region $\{f \leq \rho\}$. Therefore, (\ref{upper.bound}) is satisfied in this case.

\textit{Case 2:} Suppose that $p_0 \in \{f=\rho\}$. Since $f^2 \, \text{\rm Ric} \geq c \, g$, we have 
\[\frac{\theta \, c}{2} \leq \theta \, \rho^2 \, \text{\rm Ric}(e_1,e_1) - \frac{\theta \, c}{2} = \rho^2 \, h(e_1,e_1) \leq \rho^2 \, \sup_{\{f=\rho\}} |h|.\] 
Since $h \leq \theta \, (\text{\rm Ric} - \frac{c}{2} \, \rho^{-2} \, g)$, we conclude that 
\[h \leq C \, \rho^2 \, \Big ( \sup_{\{f=\rho\}} |h| \Big ) \, g\] 
at each point in the region $\{f \leq \rho\}$. This proves (\ref{upper.bound}). \\

\begin{lemma}
\label{max.principle.h}
Let $h$ be a solution of the Lichnerowicz-type equation 
\[\Delta_L h + \mathscr{L}_X(h) = 0\] 
on the region $\{f \leq \rho\}$. Then 
\[\sup_{\{f \leq \rho\}} f^2 \, |h| \leq B \, \rho^2 \sup_{\{f=\rho\}} |h|,\] 
where $B$ is a positive constant that does not depend on $\rho$.
\end{lemma}

\textbf{Proof.} 
As above, it suffices to show that 
\begin{equation} 
\label{upper.bound.2}
f^2 \, h \leq C \, \rho^2 \, \Big ( \sup_{\{f=\rho\}} |h| \Big ) \, g 
\end{equation}
for some uniform constant $C$. We now describe the proof of (\ref{upper.bound.2}). By Proposition \ref{asympt}, we can find a compact set $K$ such that $f \, \text{\rm Ric} < (1 - 3 \, f^{-1} \, |\nabla f|^2) \, g$ on $M \setminus K$. Let us consider the smallest real number $\theta$ with the property that $\theta \, f^{-2} \, g - h$ is positive semi-definite at each point in the region $\{f \leq \rho\}$. By definition of $\theta$, there exists a point $p_0 \in \{f \leq \rho\}$ and an orthonormal basis $\{e_1,\hdots,e_n\}$ of $T_{p_0} M$ such that 
\[\theta \, f^{-2} - h(e_1,e_1) = 0\] 
at the point $p_0$. Let us distinguish two cases: 

\textit{Case 1:} Suppose that $p_0 \in \{f < \rho\} \setminus K$. In this case, we have 
\[\theta \, \Delta(f^{-2}) - (\Delta h)(e_1,e_1) \geq 0\] 
and 
\[\theta \, \langle X,\nabla(f^{-2}) \rangle - (D_X h)(e_1,e_1) = 0\] 
at the point $p_0$. Using the identity $\Delta_L h + \mathscr{L}_X(h) = 0$, we obtain 
\begin{align*} 
0 &= (\Delta h)(e_1,e_1) + (D_X h)(e_1,e_1) + 2 \sum_{i,k=1}^n R(e_1,e_i,e_1,e_k) \, h(e_i,e_k) \\ 
&\leq \theta \, \Delta (f^{-2}) + \theta \, \langle X,\nabla(f^{-2}) \rangle + 2 \sum_{i,k=1}^n R(e_1,e_i,e_1,e_k) \, h(e_i,e_k) \\ 
&= -2\theta \, f^{-3} \, (1 - 3 \, f^{-1} \, |\nabla f|^2 - f \, \text{\rm Ric}(e_1,e_1)) \\ 
&- 2 \sum_{i,k=1}^n R(e_1,e_i,e_1,e_k) \, (\theta \, f^{-2} \, g(e_i,e_k) - h(e_i,e_k)) 
\end{align*} 
at the point $p_0$. Since $(M,g)$ has positive sectional curvature, we have 
\[\sum_{i,k=1}^n R(e_1,e_i,e_1,e_k) \, (\theta \, f^{-2} \, g(e_i,e_k) - h(e_i,e_k)) \geq 0,\] 
hence 
\[0 \leq -2\theta \, f^{-3} \, (1 - 3 \, f^{-1} \, |\nabla f|^2 - f \, \text{\rm Ric}(e_1,e_1)).\] 
On the other hand, we have $f \, \text{\rm Ric}(e_1,e_1) < 1 - 3 \, f^{-1} \, |\nabla f|^2$ since $p_0 \in M \setminus K$. Consequently, we have $\theta \leq 0$. This implies that $h \leq 0$ at each point in the region $\{f \leq \rho\}$, and (\ref{upper.bound.2}) is trivially satisfied. 

\textit{Case 2:} We next assume that $p_0 \in \{f=\rho\} \cup K$. Using Lemma \ref{prep}, we obtain  
\[\theta = f^2 \, h(e_1,e_1) \leq \sup_{\{f=\rho\} \cup K} f^2 \, |h| \leq C \, \rho^2 \, \sup_{\{f=\rho\}} |h|.\] 
Since $f^2 \, h \leq \theta \, g$, we conclude that 
\[f^2 \, h \leq C \, \rho^2 \, \Big ( \sup_{\{f=\rho\}} |h| \Big ) \, g\] 
at each point in the region $\{f \leq \rho\}$. This proves (\ref{upper.bound.2}). \\

\begin{theorem}
\label{lichnerowicz.equation}
Suppose that $h$ is a solution of the Lichnerowicz-type equation 
\[\Delta_L h + \mathscr{L}_X(h) = 0.\] 
with the property that $|h| \leq O(r^{\frac{1}{2(n-2)}-\frac{1}{2}-\varepsilon})$. Then $h = \lambda \, \text{\rm Ric}$ for some constant $\lambda \in \mathbb{R}$.
\end{theorem} 

\textbf{Proof.} Let us consider the function 
\[A(r) = \inf_{\lambda \in \mathbb{R}} \sup_{\{f=r\}} |h - \lambda \, \text{\rm Ric}|.\] 
Clearly, $A(r) \leq \sup_{\{f=r\}} |h| \leq O(r^{\frac{1}{2(n-2)}-\frac{1}{2}-\varepsilon})$. We consider two cases: 

\textit{Case 1:} Suppose that there exists a sequence of real numbers $r_m \to \infty$ such that $A(r_m) = 0$ for all $m$. In this case, we can find a sequence of real numbers $\lambda_m$ such that $h - \lambda_m \, \text{\rm Ric} = 0$ on the level surface $\{f=r_m\}$. Using Lemma \ref{max.principle.h}, we conclude that $h - \lambda_m \, \text{\rm Ric} = 0$ in the region $\{f \leq r_m\}$. Therefore, the sequence $\lambda_m$ is constant. Moreover, $h$ is a constant multiple of the Ricci tensor.

\textit{Case 2:} Suppose now that $A(r) > 0$ when $r$ is sufficiently large. Let us fix a real number $\tau \in (0,\frac{1}{2})$ such that $\tau^{-\varepsilon} > 2N \, B$, where $N$ and $B$ are the constants in Lemma \ref{spherical.harmonics.h} and Lemma \ref{max.principle.h}, respectively. Since $A(r) \leq O(r^{\frac{1}{2(n-2)}-\frac{1}{2}-\varepsilon})$, there exists a sequence of real numbers $r_m \to \infty$ such that 
\[A(r_m) \leq 2 \, \tau^{\frac{1}{2}-\frac{1}{2(n-2)}+\varepsilon} \, A(\tau r_m)\] 
for all $m$. For each $m$, we can find a real number $\lambda_m$ such that 
\[\sup_{\{f=r_m\}} |h - \lambda_m \, \text{\rm Ric}| = A(r_m).\] 
Applying Lemma \ref{max.principle.h} to the tensor 
\[\tilde{h}^{(m)} = \frac{1}{A(r_m)} \, (h - \lambda_m \, \text{\rm Ric})\] 
gives 
\begin{equation} 
\label{bound.for.h} 
\sup_{\{f=r\}} |\tilde{h}^{(m)}| \leq \frac{B \, r_m^2}{r^2} \, \sup_{\{f=r_m\}} |\tilde{h}^{(m)}| = \frac{B \, r_m^2}{r^2 \, A(r_m)} \, \sup_{\{f=r_m\}} |h - \lambda_m \, \text{\rm Ric}| = \frac{B \, r_m^2}{r^2} 
\end{equation}
for $r \leq r_m$. 

At this point, we define 
\[\hat{g}^{(m)}(t) = r_m^{-1} \, \Phi_{r_m t}^*(g)\] 
and 
\[\hat{h}^{(m)}(t) = r_m^{-1} \, \Phi_{r_m t}^*(\tilde{h}^{(m)}).\] 
The metrics $\hat{g}^{(m)}(t)$ evolve by the Ricci flow, and the tensors $\hat{h}^{(m)}(t)$ satisfy the parabolic Lichnerowicz equation 
\[\frac{\partial}{\partial t} \hat{h}^{(m)}(t) = \Delta_{L,\hat{g}^{(m)}(t)} \hat{h}^{(m)}(t).\] 
Using (\ref{bound.for.h}), we obtain 
\[\limsup_{m \to \infty} \sup_{t \in [\delta,1-\delta]} \sup_{\{r_m - \delta^{-1} \, \sqrt{r_m} \leq f \leq r_m + \delta^{-1} \, \sqrt{r_m}\}} |\hat{h}^{(m)}(t)|_{\hat{g}^{(m)}(t)} < \infty\] 
for any given $\delta \in (0,\frac{1}{2})$. 

We now pass to the limit as $m \to \infty$. Let us choose a sequence of marked points $p_m \in M$ satisfying $f(p_m) = r_m$. The manifolds $(M,\hat{g}^{(m)}(t),p_m)$ converge in the Cheeger-Gromov sense to a one-parameter family of shrinking cylinders $(S^{n-1} \times \mathbb{R},\overline{g}(t))$, $t \in (0,1)$, where $\overline{g}(t)$ is given by (\ref{shrinking.cylinders}). The vector fields $r_m^{\frac{1}{2}} \, X$ converge to the axial vector field $\frac{\partial}{\partial z}$ on $S^{n-1} \times \mathbb{R}$. Furthermore, the sequence $\hat{h}^{(m)}(t)$ converges to a one-parameter family of tensors $\overline{h}(t)$, $t \in (0,1)$, which solve the parabolic Lichnerowicz equation 
\[\frac{\partial}{\partial t} \overline{h}(t) = \Delta_{L,\overline{g}(t)} \overline{h}(t).\] 
As in \cite{Brendle2}, we can show that $\overline{h}(t)$ is invariant under translations along the axis of the cylinder. Using (\ref{bound.for.h}), we obtain 
\[|\overline{h}(t)|_{\overline{g}(t)} \leq B \, (1-t)^{-2}\] 
for all $t \in (0,\frac{1}{2}]$. Hence, Lemma \ref{spherical.harmonics.h} implies that 
\begin{equation} 
\label{symmetry.h}
\inf_{\lambda \in \mathbb{R}} \sup_{S^{n-1} \times \mathbb{R}} \big | \overline{h}(t) - \lambda \, \text{\rm Ric}_{\overline{g}(t)} \big |_{\overline{g}(t)} \leq N \, B \, (1-t)^{\frac{1}{2(n-2)}-\frac{1}{2}} 
\end{equation}
for all $t \in [\frac{1}{2},1)$. On the other hand, we have 
\begin{align*} 
&\inf_{\lambda \in \mathbb{R}} \sup_{\Phi_{r_m(\tau-1)}(\{f=\tau r_m\})} \Big | \hat{h}^{(m)}(1-\tau) - \lambda \, \text{\rm Ric}_{\hat{g}^{(m)}(1-\tau)} \Big |_{\hat{g}^{(m)}(1-\tau)} \\ 
&= \inf_{\lambda \in \mathbb{R}} \sup_{\{f=\tau r_m\}} |\tilde{h}^{(m)} - \lambda \, \text{\rm Ric}_g|_g \\ 
&= \frac{1}{A(r_m)} \, \inf_{\lambda \in \mathbb{R}} \sup_{\{f=\tau r_m\}} |h - \lambda \, \text{\rm Ric}_g|_g \\ 
&= \frac{A(\tau r_m)}{A(r_m)} \\ 
&\geq \frac{1}{2} \, \tau^{\frac{1}{2(n-2)}-\frac{1}{2}-\varepsilon}. 
\end{align*}
If we send $m \to \infty$, we obtain 
\begin{equation} 
\label{asymmetry.h}
\inf_{\lambda \in \mathbb{R}} \sup_{S^{n-1} \times \mathbb{R}} \big | \overline{h}(1-\tau) - \lambda \, \text{\rm Ric}_{\overline{g}(1-\tau)} \big |_{\overline{g}(1-\tau)} \geq \frac{1}{2} \, \tau^{\frac{1}{2(n-2)}-\frac{1}{2}-\varepsilon}. 
\end{equation} 
Since $\tau^{-\varepsilon} > 2N \, B$, the inequality (\ref{asymmetry.h}) contradicts (\ref{symmetry.h}). This completes the proof of Theorem \ref{lichnerowicz.equation}. \\

\section{Proof of Theorem \ref{main.theorem}}

\label{sym.prin}

Combining Theorems \ref{existence.of.V} and \ref{lichnerowicz.equation}, we obtain the following symmetry principle: 

\begin{theorem}
\label{symmetry.principle}
Suppose that $U$ is a vector field on $(M,g)$ such that $|\mathscr{L}_U(g)| \leq O(r^{\frac{1}{2(n-2)}-\frac{1}{2}-2\varepsilon})$ and $|\Delta U + D_X U| \leq O(r^{\frac{1}{2(n-2)}-1-2\varepsilon})$ for some small constant $\varepsilon>0$. Then there exists a vector field $\hat{U}$ on $(M,g)$ such that $\mathscr{L}_{\hat{U}}(g) = 0$, $[\hat{U},X] = 0$, $\langle \hat{U},X \rangle = 0$, and $|\hat{U} - U| \leq O(r^{\frac{1}{2(n-2)}-\varepsilon})$.
\end{theorem}

\textbf{Proof.} 
In view of Theorem \ref{existence.of.V}, the equation 
\[\Delta V + D_X V = \Delta U + D_X U\] 
has a smooth solution which satisfies the bounds $|V| \leq O(r^{\frac{1}{2(n-2)}-\varepsilon})$ and $|DV| \leq O(r^{\frac{1}{2(n-2)}-\frac{1}{2}-\varepsilon})$. Hence, the vector field $W = U-V$ satisfies $\Delta W + D_X W = 0$. Using Theorem 4.1 in \cite{Brendle2}, we conclude that the Lie derivative $h = \mathscr{L}_W(g)$ satisfies the Lichnerowicz-type equation 
\[\Delta_L h + \mathscr{L}_X(h) = 0.\] 
Since $|h| \leq O(r^{\frac{1}{2(n-2)}-\frac{1}{2}-\varepsilon})$, Theorem \ref{lichnerowicz.equation} implies that $h = \lambda \, \text{\rm Ric}$ for some constant $\lambda \in \mathbb{R}$. Consequently, the vector field $\hat{U} := U - V - \frac{1}{2} \, \lambda \, X$ must be a Killing vector field. The identities $[\hat{U},X] = 0$ and $\langle \hat{U},X \rangle = 0$ follow as in \cite{Brendle2}. \\

To complete the proof of Theorem \ref{main.theorem}, we apply Theorem \ref{symmetry.principle} to the vector fields $U_a$ constructed in Proposition \ref{rotation.vector.fields}. Consequently, there exist vector fields $\hat{U}_a$, $a \in \{1,\hdots,\frac{n(n-1)}{2}\}$, on $(M,g)$ such that $\mathscr{L}_{\hat{U}_a}(g) = 0$, $[\hat{U}_a,X] = 0$, and $\langle \hat{U}_a,X \rangle = 0$. Moreover, we have 
\[\sum_{a=1}^{\frac{n(n-1)}{2}} \hat{U}_a \otimes \hat{U}_a = r \, \Big ( \sum_{i=1}^{n-1} e_i \otimes e_i + O(r^{\frac{1}{2(n-2)}-\frac{1}{2}-\varepsilon}) \Big ),\] 
where $\{e_1,\hdots,e_{n-1}\}$ is a local orthonormal frame on the level set $\{f=r\}$. This shows that $(M,g)$ is rotationally symmetric.

\section{Proof of Theorem \ref{dim.4}}

We now describe how Theorem \ref{dim.4} follows from Theorem \ref{main.theorem}. Let $(M,g)$ be a four-dimensional steady gradient Ricci soliton which is non-flat; is $\kappa$-noncollapsed; and satisfies the pointwise pinching condition 
\[0 \leq \max \{a_3,b_3,c_3\} \leq \Lambda \, \min \{a_1+a_2,c_1+c_2\}\] 
for some constant $\Lambda \geq 1$. In particular, $(M,g)$ has nonnegative isotropic curvature. Moreover, since the sum $R + |\nabla f|^2$ is constant, the scalar curvature of $(M,g)$ is bounded from above; consequently, $(M,g)$ has bounded curvature.

We next show that $(M,g)$ has positive curvature operator. To that end, we adapt the arguments in \cite{Hamilton1} and \cite{Hamilton2}. We note that pinching estimates for ancient solutions to the Ricci flow were established in \cite{Brendle-Huisken-Sinestrari}.

\begin{lemma}
\label{pinching.1}
We have $a_3 \leq (6\Lambda^2+1) \, a_1$ and $c_3 \leq (6\Lambda^2+1) \, a_1$.
\end{lemma} 

\textbf{Proof.} 
Using the inequalities 
\[\Delta a_1 + \langle X,\nabla a_1 \rangle \leq -2a_2a_3\] 
and 
\[\Delta a_3 + \langle X,\nabla a_3 \rangle \geq -a_3^2 - 2a_1a_2 - b_3^2,\] 
we obtain 
\begin{align*} 
&\Delta ((6\Lambda^2+1) \, a_1 - a_3) + \langle X,\nabla ((6\Lambda^2+1) \, a_1 - a_3) \rangle \\ 
&\leq a_3^2+2a_1a_2+b_3^2 - (12\Lambda^2+2) \, a_2a_3 \\ 
&\leq a_3^2+b_3^2 - 12\Lambda^2 \, a_2a_3 \\ 
&\leq a_3^2+b_3^2 - 3\Lambda^2 \, (a_1+a_2)^2 \\ 
&\leq -a_3^2. 
\end{align*} 
Hence, the Omori-Yau maximum principle implies that $(6\Lambda^2+1) \, a_1 - a_3 \geq 0$. The inequality $(6\Lambda^2+1) \, c_1 - c_3 \geq 0$ follows similarly. \\

\begin{lemma}
\label{pinching.2}
We have $4b_3^2 \leq (a_1+a_2) \, (c_1+c_2)$.
\end{lemma}

\textbf{Proof.} 
Suppose that $\gamma = \sup_M \frac{2b_3}{\sqrt{(a_1+a_2) \, (c_1+c_2)}} > 1$. The function $u = \frac{1}{2} \, \sqrt{(a_1+a_2) \, (c_1+c_2)}$ satisfies 
\begin{align*} 
&\Delta u + \langle X,\nabla u \rangle \\ 
&\leq -u \, \bigg [ a_3+c_3 + \frac{a_1^2+a_2^2+b_1^2+b_2^2}{2(a_1+a_2)} + \frac{c_1^2+c_2^2+b_1^2+b_2^2}{2(c_1+c_2)} \bigg ]. 
\end{align*} 
On the other hand, we have 
\[\Delta b_3 + \langle X,\nabla b_3 \rangle \geq -b_3(a_3+c_3) - 2b_1b_2.\] 
Putting these facts together, we obtain 
\begin{align*} 
&\Delta (\gamma u-b_3) + \langle X,\nabla (\gamma u-b_3) \rangle \\ 
&\leq -\gamma u \, \bigg [ a_3+c_3 + \frac{a_1^2+a_2^2+b_1^2+b_2^2}{2(a_1+a_2)} + \frac{c_1^2+c_2^2+b_1^2+b_2^2}{2(c_1+c_2)} \bigg ] \\ 
&+ b_3 (a_3+c_3) + 2b_1 b_2 \\ 
&= -\gamma u \, \frac{(a_1-b_1)^2 + (a_2-b_2)^2 + 2a_2(b_2-b_1)}{2(a_1+a_2)} \\ 
&- \gamma u \, \frac{(c_1-b_1)^2 + (c_2-b_2)^2 + 2c_2(b_2-b_1)}{2(c_1+c_2)} \\ 
&- (\gamma u - b_3) \, (a_3+c_3 + 2b_1) - 2b_1(b_3-b_2). 
\end{align*} 
Note that $\gamma u - b_3 \geq 0$ by definition of $\gamma$. Since $\gamma > 1$, we can find a positive constant $\delta$ such that 
\begin{align*}
3\delta \, |\text{\rm Ric}|^2 
&\leq \gamma u \, \frac{(a_1-b_1)^2 + (a_2-b_2)^2 + 2a_2(b_2-b_1)}{2(a_1+a_2)} \\ 
&+ \gamma u \, \frac{(c_1-b_1)^2 + (c_2-b_2)^2 + 2c_2(b_2-b_1)}{2(c_1+c_2)} \\ 
&+ (\gamma u - b_3) \, (a_3+c_3 + 2b_1) + 2b_1(b_3-b_2). 
\end{align*} 
This implies 
\[\Delta (\gamma u-b_3) + \langle X,\nabla (\gamma u-b_3) \rangle \leq -3\delta \, |\text{\rm Ric}|^2,\] 
hence 
\[\Delta (\gamma u-b_3-\delta \, R) + \langle X,\nabla (\gamma u-b_3-\delta \, R) \rangle \leq -\delta \, |\text{\rm Ric}|^2.\] 
Using the Omori-Yau maximum principle, we conclude that $\gamma u - b_3 - \delta \, R \geq 0$. This contradicts the definition of $\gamma$. Thus, $\gamma \leq 1$, as claimed. \\

\begin{proposition}
\label{pinching.3}
We have $b_3^2 \leq a_1c_1$.
\end{proposition}

\textbf{Proof.} 
Suppose that $\gamma = \sup_M \frac{b_3}{\sqrt{a_1c_1}} > 1$. The function $v = \sqrt{a_1 c_1}$ satisfies 
\[\Delta v + \langle X,\nabla v \rangle \leq -v \, \bigg [ \frac{a_1^2+2a_2a_3+b_1^2}{2a_1} + \frac{c_1^2+2c_2c_3+b_1^2}{2c_1} \bigg ].\] 
This implies 
\begin{align*} 
&\Delta (\gamma v - b_3) + \langle X,\nabla (\gamma v - b_3) \rangle \\ 
&\leq -\gamma v \, \bigg [ \frac{a_1^2+2a_2a_3+b_1^2}{2a_1} + \frac{c_1^2+2c_2c_3+b_1^2}{2c_1} \bigg ] \\ 
&+ b_3 (a_3+c_3) + 2b_1 b_2 \\ 
&= -\gamma v \, \bigg [ \frac{(a_1-b_1)^2+2(a_2-a_1)a_3}{2a_1} + \frac{(c_1-b_1)^2+2(c_2-c_1)c_3}{2c_1} \bigg ] \\ 
&- (\gamma v - b_3) \, (a_3+c_3+2b_1) - 2b_1(b_3-b_2). 
\end{align*} 
Note that $\gamma v - b_3 \geq 0$ by definition of $\gamma$. Using Lemma \ref{pinching.2} and the inequality $\gamma>1$, we obtain an estimate of the form 
\begin{align*}
3\delta \, |\text{\rm Ric}|^2 
&\leq \gamma v \, \bigg [ \frac{(a_1-b_1)^2+2(a_2-a_1)a_3}{2a_1} + \frac{(c_1-b_1)^2+2(c_2-c_1)c_3}{2c_1} \bigg ] \\ 
&+ (\gamma v - b_3) \, (a_3+c_3+2b_1) + 2b_1(b_3-b_2) 
\end{align*} 
for some positive constant $\delta$. From this, we deduce that 
\[\Delta (\gamma v-b_3) + \langle X,\nabla (\gamma v-b_3) \rangle \leq -3\delta \, |\text{\rm Ric}|^2,\] 
hence 
\[\Delta (\gamma v-b_3-\delta \, R) + \langle X,\nabla (\gamma v-b_3-\delta \, R) \rangle \leq -\delta \, |\text{\rm Ric}|^2.\] 
As above, the Omori-Yau maximum principle implies that $\gamma v - b_3 - \delta \, R \geq 0$. This contradicts the definition of $\gamma$. Consequently, $\gamma \leq 1$, which proves the assertion. \\

\begin{corollary}
\label{pco}
The manifold $(M,g)$ has positive curvature operator.
\end{corollary}

\textbf{Proof.} 
The inequality $b_3^2 \leq a_1c_1$ implies that $(M,g)$ has nonnegative curvature operator. If $(M,g)$ has generic holonomy group, then the strict maximum principle (cf. \cite{Hamilton1}) implies that $(M,g)$ has positive curvature operator. On the other hand, if $(M,g)$ has non-generic holonomy group, then $(M,g)$ locally splits as a product. In this case, we can deduce from Proposition \ref{pinching.3} that $(M,g)$ is isometric to a cylinder. This contradicts the fact that $(M,g)$ is a steady soliton. \\

Note that $(M,g)$ satisfies restricted isotropic curvature pinching condition in \cite{Chen-Zhu}. Using the compactness theorem for ancient $\kappa$-solutions in \cite{Chen-Zhu}, we obtain:

\begin{proposition}[B.L.~Chen and X.P.~Zhu \cite{Chen-Zhu}]
\label{basics}
Let $p_m$ be a sequence of points going to infinity. Then $|\langle X,\nabla R \rangle| \leq O(1) \, R^2$ at the point $p_m$. Moreover, if $d(p_0,p_m)^2 \, R(p_m) \to \infty$, then we have $|\nabla R| \leq o(1) \, R^{\frac{3}{2}}$ and $|\langle X,\nabla R \rangle + \frac{2}{3} \, R^2| \leq o(1) \, R^2$ at the point $p_m$.
\end{proposition}

\textbf{Proof.} 
This is a consequence of Lemma 3.1, Proposition 3.3, and Corollary 3.7 in \cite{Chen-Zhu}. \\

Using Proposition \ref{basics}, it is not difficult to show that $R \to 0$ at infinity. Consequently, there exists a unique point $p_0 \in M$ where the scalar curvature attains its maximum. The point $p_0$ must be the a critical point of the function $f$. Since $f$ is strictly convex, we conclude that $f$ grows linearly near infinity. If we integrate the inequality $|\langle X,\nabla R \rangle| \leq O(1) \, R^2$ along integral curves of $X$, we obtain $R \geq \frac{\Lambda_1}{d(p_0,p)}$ outside a compact set, where $\Lambda_1$ is a positive constant. Hence, Proposition \ref{basics} gives $|\langle X,\nabla R \rangle + \frac{2}{3} \, R^2| \leq o(1) \, R^2$. Integrating this inequality along integral curves of $X$, we obtain $R \leq \frac{\Lambda_2}{d(p_0,p)}$ outside a compact set. Using Lemma 3.1 in \cite{Chen-Zhu} again, we conclude that $(M,g)$ is asymptotically cylindrical. Hence, $(M,g)$ must be rotationally symmetric by Theorem \ref{main.theorem}.

\appendix

\section{The eigenvalues of some elliptic operators on $S^{n-1}$}

In this section, we analyze the eigenvalues of certain elliptic operators on $S^{n-1}$. In the following, $g_{S^{n-1}}$ will denote the standard metric on $S^{n-1}$ with constant sectional curvature $1$. 

\begin{proposition} 
\label{laplacian.on.one.forms}
Let $\sigma$ be a one-form on $S^{n-1}$ satisfying 
\[\Delta_{S^{n-1}} \sigma + \mu \, \sigma = 0,\] 
where $\Delta_{S^{n-1}}$ denotes the rough Laplacian and $\mu \in (-\infty,1)$ is a constant. Then $\sigma = 0$. 
\end{proposition} 

\textbf{Proof.} 
For any smooth function $u$, we have 
\begin{align*} 
\int_{S^{n-1}} u \, \Delta_{S^{n-1}}(d^* \sigma) 
&= \int_{S^{n-1}} \langle d(\Delta_{S^{n-1}} u),\sigma \rangle \\ 
&= \int_{S^{n-1}} \langle \Delta_{S^{n-1}} (du),\sigma \rangle - (n-2) \int_{S^{n-1}} \langle du,\sigma \rangle \\ 
&= -(n-2+\mu) \int_{S^{n-1}} \langle du,\sigma \rangle \\ 
&= -(n-2+\mu) \int_{S^{n-1}} u \, d^* \sigma. 
\end{align*} 
Since $u$ is arbitrary, we conclude that 
\[\Delta_{S^{n-1}} (d^* \sigma) + (n-2+\mu) \, d^* \sigma = 0.\] 
Since $n-2+\mu < n-1$, it follows that $d^* \sigma$ is constant. Consequently, $d^* \sigma = 0$ by the divergence theorem. 

We next consider the tensor $S_{ik} = D_i \sigma_k + D_k \sigma_i$. Then 
\[(n-2-\mu) \, \sigma_i = \Delta_{S^{n-1}} \sigma_i + (n-2) \, \sigma_i = D^k S_{ik} - \frac{1}{2} \, D_i(\text{\rm tr} \, S).\] 
Using the identity $d^* \sigma = 0$, we obtain 
\[(n-2-\mu) \int_{S^{n-1}} |\sigma|^2 = \int_{S^{n-1}} \Big ( D^k S_{ik} - \frac{1}{2} \, D_i(\text{\rm tr} \, S) \Big ) \, \sigma^i = -\frac{1}{2} \int_{S^{n-1}} |S|^2.\] 
Since $n-2-\mu > 0$, we conclude that $\sigma = 0$, as claimed. \\

\begin{proposition} 
\label{laplacian.on.tensors}
Let $\chi$ be a symmetric $(0,2)$-tensor on $S^{n-1}$ satisfying 
\[\Delta_{S^{n-1}} \chi - 2(n-1) \, \overset{\text{\rm o}}{\chi} + \mu \, \chi = 0,\] 
where $\overset{\text{\rm o}}{\chi}$ denotes the trace-free part of $\chi$ and $\mu \in (-\infty,n-1)$ is a constant. Then $\chi$ is a constant multiple of $g_{S^{n-1}}$.
\end{proposition} 

\textbf{Proof.} 
The trace of $\chi$ satisfies 
\[\Delta_{S^{n-1}}(\text{\rm tr} \, \chi) + \mu \, (\text{\rm tr} \, \chi) = 0.\] 
Since $\mu < n-1$, we conclude that $\text{\rm tr} \, \chi$ is constant. Moreover, the trace-free part of $\chi$ satisfies 
\[\Delta_{S^{n-1}} \overset{\text{\rm o}}{\chi} + (\mu - 2(n-1)) \, \overset{\text{\rm o}}{\chi} = 0.\] 
Since $\mu-2(n-1) < 0$, it follows that $\overset{\text{\rm o}}{\chi} = 0$. Putting these facts together, the assertion follows.

\end{document}